# DEA-based benchmarking for performance evaluation in pay-for-performance incentive plans

*Wade D. Cook[1], Nuria Ramón[2], José L. Ruiz[2], Inmaculada Sirvent[2] and Joe Zhu[3]*

[1] *Schulich School of Business, York University, Toronto, Ontario M3J 1P3, Canada*
[2] *Centro de Investigación Operativa, Universidad Miguel Hernández, Avd. de la Universidad, s/n, 03202-Elche (Alicante), Spain*
[3] *School of Business, Worcester Polytechnic Institute, Worcester, MA 01609, USA*

*February, 2018*

## Abstract

Incentive plans involve payments for performance relative to some set of goals. In this paper, we extend Data Envelopment Analysis (DEA) to the evaluation of performance in the specific context of pay-for-performance incentive plans. The approach proposed ensures that the evaluation of performance of decision making units (DMUs) that follow the implementation of incentive plans, is made in terms of targets that are attainable, as well as representing best practices. A model is developed that adjusts the benchmarking to the goals through the corresponding payment of incentives, thus DEA targets are established taking into consideration the improvement strategies that were set out in the incentive plans. To illustrate, we examine an application concerned with the financing of public Spanish universities.

## 1. Introduction

This paper extends DEA to incorporate management goals in the evaluation of the performance of DMUs. Specifically, the approach proposed is developed within the context of pay-for-performance incentive plans, so that the ultimate aim is to provide a benchmarking framework based on appropriate targets for the monitoring of improvements.

In planning improvements, organizations often implement incentive plans, which typically pay for performance relative to some goals (some studies have shown that performance improves as a result of the implementation of these plans; see, e.g., Banker et al. (1996, 2000)). Once the goals have been set, the payment of incentives is carried out in a subsequent period of performance, as a result of an evaluation of the level of achievement of goals. Payment can, for example, be tied to performance by means of a linear function connecting incentive rates and the degree of achievement of goals. The financing of public Spanish universities can be seen as an example of a situation wherein some organizations follow the implementation of pay-for-performance incentive plans. This is due to the fact that part of their financing is linked to the achievement of goals relating to the different areas of performance (teaching, research, knowledge transfer,…). Incentive plans of these higher education institutions are designed through programme-contracts, which formalize an agreement between regional governments and university managers, in which goal levels are set for a number of index indicators. The payment incentive is made in the next budgetary year on the basis of an evaluation of performance relative to the goals.

Implementing pay-for-performance incentive plans involves (1) a first step of planning, in which goals are set, and (2) a subsequent step of monitoring and control, in which performance is evaluated and the payment of incentives is established. This paper focuses exclusively on the latter step involving the evaluation of performance, without regard to the way goals were laid out. Obviously, such evaluation is, however, determined by the goal levels previously set. In this respect, it should be highlighted that, on one hand, goals might have been set without having any evidence that they will be achievable at the moment of conducting performance evaluation, and on the other, that in participative goal setting, managers whose own performance is to be evaluated have the opportunity to influence the goals (downward) in order to maximize the payout (see Anderson et al. (2010) for discussions). Taking into account these considerations, we develop an approach for evaluating performance in the context of monitoring incentive plans, by using targets that are attainable, and as well represent best practices. Furthermore, we seek targets that respect, as much as possible, the payment of incentives relative to the goals.

In order to accomplish the above, we propose the use of Data Envelopment Analysis (DEA) (Charnes et al., 1978), which is a tool directed to evaluating past performances as part of the control function of management (as stated in Cooper (2005)). DEA evaluates performance of decision making units (DMUs) from a perspective of benchmarking, through the setting of targets on a best-practice frontier, thereby making it possible to evaluate them against their peers[1]. See also Thanassoulis et al. (2008) and Cook et al. (2014) for discussions on the issue of benchmarking in DEA. In the specific situation of having a set of DMUs that follow the implementation of incentive plans, DEA can provide an alternative evaluation of performances in time period t, based on targets that, in addition to being attainable and efficient, have been established taking into account the goals which were set in period t-1. Specifically, the model we develop here finds DEA targets by minimizing the differences between the payments of incentives relative to DEA targets and those relative to the goals. Relating DEA targets and goals seeks to consider the policy of improvements that was pursued through the setting of goals (note that DMUs might have oriented their activities towards the achievement of the goals, in particular because they have a monetary incentive). Thus, as in a conventional DEA benchmarking analysis, we eventually set targets, and the comparisons between them and actual data allow us to evaluate performances (which includes analyzing individual strengths and weaknesses and establishing directions for improvement) as part of the monitoring process of the pay-for-performance incentive plans. See Stewart (2010) which also proposes an approach for setting targets within a general framework that incorporates long term goals to the DEA models. As with our approach, that author looks for realistically achievable targets, which are on the efficient frontier, and are found taking into account the goals. Nevertheless, his approach is aimed at planning, while ours is intended mainly for monitoring and control.

[1] In this sense, DEA is a tool that can be used not only in management control but also it provides certain degrees of support to management planning, insofar as it can provide suggestions as to where and by how much an inefficient DMU should be improved in order to achieve full efficiency in comparison to its peers (as discussed in Yang et al. (2009)).

The approach proposed herein is developed in the context of performance evaluation of DMUs in terms of a number of index indicators. Index indicators are frequently used in the assessment of health services, university performance, wealth or human development in hospitals, higher education institutions, and countries. Following Liu et al. (2011), the models for the benchmarking in pay-for-performance incentive plans we formulate are DEA-type models without inputs, which can be used in the evaluation of performance without regard to any production process. It is also highlighted that the models are developed within a general framework that would allow us to deal with the specific situation of having DMUs organized into groups of units that experience similar circumstances (similar objectives, policies,…). In those cases, we do not allow for individual circumstances within the set of DMUs of a given group, in line with the within-group common benchmarking approach in Cook et al. (2017). This paper is partly motivated by the situation already mentioned of the financing of public universities in Spain, wherein many of the competences regarding university policies are transferred to the regional governments. As an illustration of the approach proposed, we examine an example on university performance in which those institutions in the same region are treated uniformly in the sense that targets result from their benchmarking against a common reference set of universities.

The paper unfolds as follows: in section 2 we explain how the payment of incentives is tied to performance in our approach. In section 3 we formulate the benchmarking model that provides the targets to be used for the evaluation of performance in the context of incentive plans. In particular, the developments that lead to a linear problem that allows us to find the optimal solution of the benchmarking model, are included in a subsection. Section 4 presents an empirical illustration of the proposed approach. Conclusions follow in section 5.

## 2. The payment of incentives

Throughout the paper, we consider that we have n comparable DMUs that follow the implementation of pay-for-performance incentive plans. We suppose that the performance of these DMUs is evaluated in terms of s index indicators $y_1,\ldots,y_s$. The DMUs can thus be described by means of the output vectors $Y_j$, where $Y_j = (y_{1j},...,y_{sj})' \geq 0$, $Y_j \neq 0$, $j = 1,...,n$, record the values of these indicators for $DMU_j$, j=1,…,n, corresponding to the performances in time period t, that is, the period in which incentives will be paid. It is assumed that larger values of these indicators represent better performances. Payment of incentives is usually tied to the degree of achievement of certain goals relative to the indicators. Denote by $\hat{y}_{1j},...,\hat{y}_{sj}$, any goals for $DMU_j$ corresponding to $y_1,\ldots,y_s$, respectively, which were set in period t-1.

It can be said that $DMU_j$ has achieved the goal set for indicator $y_r$ if $\hat{y}_{rj} \leq y_{rj}$. Otherwise, the goal would not have been reached. Nevertheless, with incentive plans, it is usually contemplated that there is the possibility of paying a certain amount of money in case of $\hat{y}_{rj} > y_{rj}$, provided that $y_{rj}$ is not too different from $\hat{y}_{rj}$. In order to address this situation, we allow for deviations $s_{rj} = \hat{y}_{rj} - y_{rj}$, r=1,…,s, j=1,…,n, and set a limit, $d_{rj}^M$, for each such deviation, in the sense that no incentives will be paid for values of $s_{rj}$ larger than or equal to $d_{rj}^M$. Formally, we can define a function $a_{rj}$, representing the degree of achievement of $DMU_j$ in goal $\hat{y}_{rj}$, in terms of the deviations, in the following manner: $a_{rj}(s_{rj}) = 1$, if $s_{rj} \leq 0$, $a_{rj}(s_{rj}) = 0$, if $s_{rj} \geq d_{rj}^M$, and $a_{rj}$ is defined as a decreasing monotonic function of $s_{rj}$ for values of the deviations in between 0 and $d_{rj}^M$. Here, the degree of achievement of goals is defined as a linear function of the deviations by setting $d_{rj}^M = y_{rj}$. That is,

$$a_{rj}(s_{rj}) = \begin{cases} 1, & s_{rj} \leq 0 \\ 1 - \dfrac{s_{rj}}{y_{rj}}, & 0 \leq s_{rj} \leq y_{rj} \\ 0, & s_{rj} \geq y_{rj} \end{cases}, \quad r=1,\ldots,s,\ j=1,\ldots,n \tag{1}$$

Thus, it is considered that the goal is not achieved if the room for improvement in the corresponding indicator is more than 100% (a different limit can be set if appropriate). In case of $0 \le s_{rj} < y_{rj}$, for example if $s_{rj}/y_{rj} = 0.2$, the room for improvement of DMU$_j$ in $y_r$ is of 20% and, therefore, we can say that its degree of achievement in that indicator is 80%, that is, $a_{rj}(s_{rj}) = 0.8$.

As for the payment itself, typically incentive rates vary linearly with the degree of achievement of goals, with a floor of zero (for cases of very poor performance; in this paper, when $s_{rj} \ge y_{rj}$) and a ceiling which can be the initially available amount of incentives (in the case of Spanish universities that amount is in some regions 10% of base financing). Formally, the payment of incentives corresponding to $y_r$ can be established through a function of the degree of achievement $a \in [0,1]$ defined as $f_{rj}(a) = Q_j \omega_r a$, where $Q_j$ are the initially available incentives for DMU$_j$ and $\omega_r$, $0 \le \omega_r \le 1$, j=1,…,n, $\sum_{r=1}^{s} \omega_r = 1$, are weights which represent the relative importance that is attached to each of the indicators.

Therefore, in terms of the deviations, the payment of incentives for DMU$_j$ corresponding to $y_r$ can be defined as $p_{rj}(s_{rj}) = f_{rj}(a_{rj}(s_{rj}))$, that is

$$p_{rj}(s_{rj}) = \begin{cases} Q_j \omega_r, & s_{rj} \le 0 \\ Q_j \omega_r \left(1 - \dfrac{s_{rj}}{y_{rj}}\right), & 0 \le s_{rj} \le y_{rj} \\ 0, & s_{rj} \ge y_{rj} \end{cases}, \quad r=1,\ldots,s,\ j=1,\ldots,n \qquad (2)$$

Graphically,

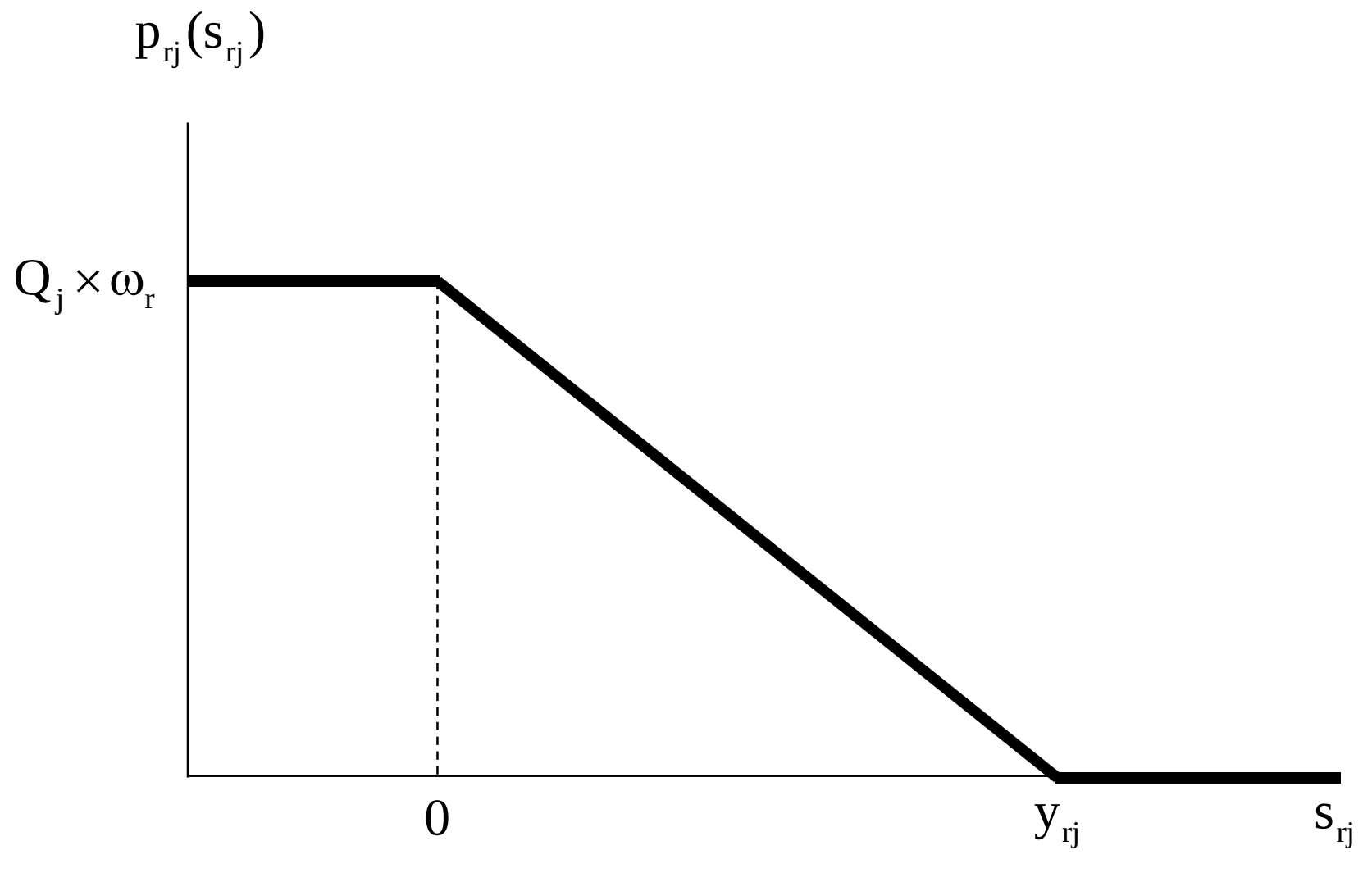

Figure 1. Payment of incentives

Eventually, the total payment of incentives for DMU$_j$ is given by $P_j(S_j) = \sum_{r=1}^{s} p_{rj}(s_{rj})$, where $S_j = (s_{1j},...,s_{sj})'$, j=1,…,n.

Remark 1. Other forms of payment that contemplate, for example, payments for overachievement of goals and penalties could be considered, by defining $p_{rj}(s_{rj})$ appropriately through the specification of new values for the limits "0" and $d_{rj}^{M}$, and/or through the addition of more ranges for $s_{rj}$ associated with the new scenarios of payment. In all cases, we should take into account that payments cannot exceed the amount of incentives that is initially available.

## 3. Goal-adjusted benchmarking in pay-for-performance incentive plans using DEA

The monitoring of improvements of a given DMU$_j$ in the context of the implementation of an incentive plan, involves the evaluation of its performance in time period t relative to the goals set in period t-1, $\hat{y}_{1j}^{g},...,\hat{y}_{sj}^{g}$. As a result, the payments are made. It is emphasized in this paper that such evaluation of performance should be made in terms of targets that are attainable and represent best practices. In addition, there is a need to relate targets and goals, in order to consider the policy of improvements that was pursued through the setting of goals. For these reasons, the approach developed here proposes adjusting the benchmarking to the goals so that targets lie on the best-practice frontier of the DEA technology associated with the actual performances of period t, while at the same time respecting as much as possible the payments of incentives relative to the goals (note that DMUs might argue that they have oriented their activities towards the achievement of the goals). In order to do so, we seek to formulate a benchmarking model that allows us to find the DEA targets, $\hat{y}_{1j}^{DEA},...,\hat{y}_{sj}^{DEA}$, which result from minimizing, for each indicator, the differences between the incentive payments relative to DEA targets, $p_{rj}(s_{rj}^{DEA})$, where $s_{rj}^{DEA} = \hat{y}_{rj}^{DEA} - y_{rj}$, and the payments relative to the goals, $p_{rj}(s_{rj}^{g})$, where $s_{rj}^{g} = \hat{y}_{rj}^{g} - y_{rj}$.

### 3.1. Model formulation

For the development to follow, the DMUs are assumed to be organized into G groups, $J_1,\ldots,J_G$, where units within each group experience similar circumstances. This means that the DMUs are assumed to be homogeneous (that is, comparable in terms of the indicators selected for the analysis), while those within the groups share objectives, policies, etc., which may differ across groups. In particular, managers may establish those common objectives and policies through setting similar goals for the DMUs within the groups. In such a situation, it may be considered appropriate to set targets treating those DMUs uniformly in terms of benchmark selection (that is, not allowing for individual circumstances). The approach proposed here provides a within-group common benchmarking framework for performance evaluation through which the DMUs within the groups are benchmarked against the same reference set, thus being treated uniformly in the identification of best practices and the setting of targets. This can be seen as a general approach, because the method developed is based on the DEA frontier determined by the entire set of DMUs. If the DMUs are not classified into groups, and we wish to allow for individual circumstances, then it can be assumed that we have as many groups as DMUs, $J_1,\ldots, J_n$, where $J_j=\{DMU_j\}$, $j=1,\ldots,n$ (as in the standard DEA). If, on the contrary, the DMUs experience similar circumstances, and we want to treat all of them uniformly in the analysis, then we can suppose that we have one group only which consists of all the DMUs, $J=\{DMU_1,\ldots,DMU_n\}$.

Since we deal with DMUs that are evaluated in terms of several index indicators, we consider for the purpose of benchmarking, the following “attainable set” proposed in Liu et al. (2011) $AS=\left\{Y\geq 0 \middle/ Y\leq\sum_{j=1}^{n}\lambda_j Y_j, \sum_{j=1}^{n}\lambda_j=1, \lambda_j\geq 0\right\}$, which is a bounded closed convex set. Liu et al. (2011) develop a theoretical framework for so-called DEA models without explicit inputs (WEI), which can be used in the evaluation of performance without regard to any production process. The benchmarking models we formulate are aimed at setting targets that result from projections of the DMUs onto the Pareto frontier of the attainable set AS. Specifically, the DMUs of a given group will all be projected onto the same facet F of the frontier of AS, because they experience similar circumstances. That is, they will be benchmarked against the same reference set of DMUs, which would be those that span that facet. Note in any case that some flexibility within the groups is permitted for target setting, because the DMUs in the groups can be projected onto different points of the facet of the frontier spanned by the common reference set (as in Cook and Zhu (2007) and Cook et al. (2017)). As said before, the choice of such facet/reference sets is made by minimizing globally the differences between the payments of incentives relative to the DEA targets and those of the goals. Summing up, DEA targets for the DMUs in a given group $J_g$ can be obtained as follows

$$\hat{Y}_j^{DEA} = Y_j + S_j^* = \left(\hat{y}_{1j}^{DEA} = y_{1j} + s_{1j}^*, ..., \hat{y}_{sj}^{DEA} = y_{sj} + s_{sj}^*\right),\ r=1,\ldots,s,\ j \in J_g \tag{3}$$

where $S_j^* = \left(s_{1j}^*, ..., s_{sj}^*\right)$, $j \in J_g$, are the optimal solutions of the model

$$\begin{array}{ll} \text{Min} & \displaystyle\sum_{j \in J_g} \sum_{r=1}^{s} \frac{\left| p_{rj}\left(s_{rj}\right) - p_{rj}\left(s_{rj}^g\right) \right|}{Q_j \omega_r} \\ \text{s.t.:} & \\ & \hat{Y}_j = Y_j + S_j = \left(\hat{y}_{1j} = y_{1j} + s_{1j}, ..., \hat{y}_{rj} = y_{rj} + s_{rj}\right) \in F \quad j \in J_g \\ & s_{rj} \text{ free} \qquad \forall r, j \end{array} \tag{4}$$

Here, F denotes a facet of the Pareto frontier of AS (see Bessent et al. (1988) and Olesen and Petersen (1996), which are two relevant papers dealing with facets in DEA). Model (4) must be solved G times, one for each group of DMUs $J_g$, g=1,…,G.

Note that in (4) the differences between the payments of incentives relative to the DEA targets and those of the goals are normalized by the amount of incentives corresponding to each index indicator initially available for each $DMU_j$, in order to avoid a possible effect of size of the DMUs.

Solving (4) will provide us with the following information:

- The targets to be used in the evaluation of performance of the DMUs. This is the key issue in the step of monitoring and control of the incentive plans, because they establish the potential in terms of which the degree of improvements of each unit is evaluated. The achievement of improvements is actually the ultimate end of the incentive plans. Regardless of the way the benchmarking is carried out, we eventually have actual performances and targets for each indicator. So, as in a conventional DEA benchmarking analysis, comparisons between these two quantities allow us to identify individual strengths and weaknesses, and define directions for further improvement.

  Moreover, comparisons between goal bundles and the bundles of DEA targets show how goals have been adjusted for the evaluation (if such adjustments actually occur).
- The payment of incentives, $P_j\left(S_j^*\right)$, for every $DMU_j$, relative to targets that are attainable and represent best practices. Obviously, large differences between $P_j\left(S_j^*\right)$ and $Q_j$ indicate that there is still room for improvement for $DMU_j$. Specifically, low incentive rates $p_{rj}\left(s_{rj}^*\right)/\left(Q_j\omega_r\right)$, r=1,…,s, show individual areas for improvement. In practice, it may happen that $P_j\left(S_j^*\right)$ and $P_j\left(S_j^g\right)$ (the incentives resulting from the evaluation relative to the goals) are similar; in fact, model (4) seeks to minimize the differences between both payments. This may serve as validation of the payments based on the goals, because in that case $DMU_j$ is given an amount of incentives that is very similar to that which would result from an appropriate target setting. However, it should be noted that, in spite of having similar payments, there might be differences between DEA targets and goals, which is important for the evaluation of the degree of improvements.

Some results are needed in order to develop an operative formulation of (4) that gives rise to the targets that are wanted. However, before presenting the technical developments in the next subsection, we explain the idea behind the approach proposed through a small example.

### 3.2. Numerical example

Consider the DMUs in Table 1, which are evaluated in terms of two indicators. For each unit, this table records the goals set in period of time t-1, the initially available amount of incentives (we suppose that half of these incentives correspond to each indicator, that is, $\omega_r = 0.5$, r=1,2, in all cases) and the actual values of the indicators in period of time t.

Table 1. Data numerical example

Table 2 displays the results provided by model (4). For each DMU, in the first row we have the actual data and the initially available amount of incentives; in the second row, the DEA targets and the payments corresponding to the evaluation made in terms of such targets; and in the third row, the goals set in t-1 and the amount of incentive corresponding to the evaluation relative to the goals. For the sake of simplicity, the setting of targets has been done separately for each unit, that is, no classification of the DMUs into groups has been considered.

Table 2. Results numerical example

Figure 2 depicts the situation graphically for DMUs D, E and F. For DMU D, the goals set in t-1 (the point $g_D$) are outside the AS corresponding to actual performances in time period t, so they can be deemed as unattainable goals. We can see how model (4) adjusts the benchmarking to these goals so that the DEA targets are the coordinates of point D', which lies on the best-practice frontier of AS. If D were evaluated with respect to the goals, it would achieve 9.17€ (of the 20€ initially available): we have a deviation from the goal of 1 unity in y1; since the actual y1 is 3, the room for improvement in this indicator is 33.3%, so it is paid 66.6% of the incentives, that is, 6.66€; following a similar reasoning, there is a room for improvement of 75% in y2, so it is paid 25% of the incentives, that is, 2.5€. Since the goals are unattainable, model (4) sets alternatively targets through D', by minimizing the differences between the payments relative to the goals and those relative to the targets (by indicators). If D were evaluated with respect to the DEA targets, which are less demanding, then it would be paid 12.17€ (following a similar reasoning to that used with the goals), 3€ more than in the evaluation relative to the goals. The evaluation of D with respect to the DEA targets reveals therefore that there is still room for improvement in both indicators for that unit (by 33.3% for y1 and 45% for y2).

For DMU E, the goals (point $g_E$) are inside AS, so they can be deemed as unambitious goals. The DEA targets on the Pareto frontier of AS (the point E') reveals that E should improve by 50% both in y1 and y2. In this case, the setting of targets results from a reduction in the payment of incentives of 4€, as DEA targets are more demanding than the goals. As for DMU F, its goals are on the Pareto frontier of AS, so DEA targets set by model (4) coincide with the goals.

Figure 2

Figure 3 depicts the situation of the DEA efficient DMUs. For DMU A, the DEA targets are adjusted to the goals, which are outside AS, through point A' on the Pareto frontier. Thus, in spite of being DEA efficient, when goals are considered, the evaluation shows that there is room for improvement for DMU A in y1 (by 100%), albeit the target has been outperformed in y2. In this case, the amount of incentives to be paid coincides with the payment relative to the goals, 12.5€ (no incentives correspond to the performance in y1 because the room for improvement is 100%, while for y2 the total available incentives are to be paid because the goal is achieved and the target is outperformed). Something similar occurs for DMU C. As for DMU B, whose goals are inside AS, DEA targets provided by model (4) coincide with actual performances, thus getting the total of available incentives.

Figure 3

### 3.3. Technical developments

In this subsection we develop an operative formulation of model (4) that allows us to find its optimal solution. In order to do so, we first reformulate (4) as follows

$$
\begin{array}{llll}
\text{Min} & \displaystyle\sum_{j\in J_g}\sum_{r=1}^{s}\frac{\left|p_{rj}\left(s_{rj}\right)-p_{rj}\left(s_{rj}^{g}\right)\right|}{Q_j\omega_r} & & \\
\text{s.t.:} & & & \\
& \displaystyle\sum_{k\in E}\lambda_{kj}y_{rk}=y_{rj}+s_{rj} & j\in J_g \quad r=1,\dots,s & (5.1)\\
& \displaystyle\sum_{k\in E}\lambda_{kj}=1 & j\in J_g & (5.2)\\
& u'Y_k+u_0+d_k=0 & k\in E & (5.3)\\
& u_r y_{rj}\ge 1 & j\in J_g \quad r=1,\dots,s & (5.4)\\
& d_k\le M b_k & k\in E & (5.5)\\
& \displaystyle\sum_{j\in J_g}\lambda_{kj}\le M(1-b_k) & k\in E & (5.6)\\
& d_k\ge 0, b_k\in\{0,1\} & k\in E & \\
& \lambda_{kj}\ge 0 & k\in E \quad j\in J_g & \\
& u_0, s_{rj} \ \text{free} & j\in J_g \quad r=1,\dots,s &
\end{array}
\qquad (5)
$$

where E is the set of extreme efficient DMUs (in the Pareto sense) of AS, and M is a large positive quantity.

The constraints (5.1)-(5.2) guarantee that the projections $\hat{Y}_j$, $j\in J_g$, belong to AS. With (5.3)-(5.4) we allow for supporting hyperplanes that contain facets of the Pareto frontier of AS (note that their coefficients are restricted to be strictly positive). (5.5)-(5.6) are the key restrictions that link the two previous groups of constraints. By virtue of these, model (5) ensures targets on the Pareto frontier of AS. Note that, if $\sum_{j\in J_g}\lambda_{kj}>0$ then $b_k=0$ and, consequently, $d_k=0$. This means that if $DMU_k$ in E participates actively as a referent in the benchmarking of some DMU$_j$, $j\in J_g$, then it necessarily belongs to $u'Y+u_0=0$. That is, the $DMU_k$'s $\in E$ that participate in the benchmarking of the DMUs of a given group $J_g$ are all on a same facet of the Pareto frontier of AS, because these $DMU_k$'s belong all to a common supporting hyperplane of AS, $u'Y+u_0=0$, whose coefficients are non-zero. Therefore, solving (5) allows us to identify a common reference set of efficient DMUs, $RS_g=\left\{DMU_\ell / \lambda_{\ell j}^{*}>0, \text{ for some } j\in J_g\right\}$, for the benchmarking of the DMUs in $J_g$, so that targets result from their projections onto the facet of the Pareto frontier of AS, spanned by the DMUs in $RS_g$. Specifically, these targets, which are actually the coordinates of the projections, can be found by using the optimal solutions of (5) as follows: $\hat{Y}_j^* = \sum_{\ell \in RS_g} \lambda_{\ell j}^* Y_\ell = \left(\hat{y}_{1j}^* = y_{1j} + s_{1j}^*, ..., \hat{y}_{sj}^* = y_{sj} + s_{sj}^*\right),\ j \in J_g$.

Remark 2. Constraints (5.5)-(5.6) include the classical big M and binary variables. Nevertheless, (5) can be solved in practice by reformulating these constraints using Special Ordered Sets (SOS) (Beale and Tomlin, 1970), which avoid the need to specify M. SOS Type 1 is a set of variables where at most one variable may be nonzero. Therefore, if we remove (5.5)-(5.6) from the formulation and define instead a SOS Type 1, for each pair of variables $\{\lambda_k, d_k\},\ k \in E$, where $\lambda_k = \sum_{j \in J_g} \lambda_{kj}$, then it is ensured that $\sum_{j \in J_g} \lambda_{kj}$ and $d_k$ cannot be simultaneously positive for $DMU_k$'s, $k \in E$. CPLEX Optimizer (and also LINGO) can solve LP problems with SOS. SOS variables have already been used for solving models like (5) in Ruiz and Sirvent (2016), Aparicio et al. (2017) and Cook et al. (2017).

The following proposition provides a linear expression of the payment of incentives $p_{rj}(s_{rj})$ as defined in (2)

**Proposition 1.** For every $DMU_j$, j=1,…n, and each indicator $y_r$, r=1,…,s, the following statement holds

$$p_{rj}(s_{rj}) = Q_j \omega_r I_{rj}^1 + Q_j \omega_r \left(1 - \frac{s_{rj}}{y_{rj}}\right) I_{rj}^2, \quad s_{rj} \in R \tag{6}$$

where $y_{rj} I_{rj}^3 - M_1 I_{rj}^1 \le s_{rj} \le y_{rj} I_{rj}^2 + M_2 I_{rj}^3$, $I_{rj}^1$, $I_{rj}^2$ and $I_{rj}^3$ being binary variables such that $I_{rj}^1 + I_{rj}^2 + I_{rj}^3 = 1$, and $M_1$ and $M_2$ are two big positive quantities.

Proof. Suppose that $p_{rj}(s_{rj})$ is defined as in (2). Let $s_{rj} < 0$. Then, according to Proposition 1, this can only happen if $I_{rj}^2 = I_{rj}^3 = 0$ and $I_{rj}^1 = 1$. In that case, (6) leads to $p_{rj}(s_{rj}) = Q_j \omega_r$, like in (2). Likewise, if $s_{rj} = 0$, we have either $I_{rj}^1 = 1$ and $I_{rj}^2 = I_{rj}^3 = 0$ or $I_{rj}^2 = 1$ and $I_{rj}^1 = I_{rj}^3 = 0$. In both cases (6) yields $p_{rj}(s_{rj}) = Q_j\omega_r$, like in (2). If $0 < s_{rj} < y_{rj}$, then $I_{rj}^2 = 1$ and $I_{rj}^1 = I_{rj}^3 = 0$, so $p_{rj}(s_{rj}) = Q_j\omega_r\left(1 - \frac{s_{rj}}{y_{rj}}\right)$, like in (2). If $s_{rj} = y_{rj}$, then either $I_{rj}^2 = 1$ and $I_{rj}^1 = I_{rj}^3 = 0$ or $I_{rj}^3 = 1$ and $I_{rj}^1 = I_{rj}^2 = 0$, so $p_{rj}(s_{rj}) = 0$, like in (2). And finally, if $s_{rj} > y_{rj}$, then $I_{rj}^3 = 1$ and $I_{rj}^1 = I_{rj}^2 = 0$, so $p_{rj}(s_{rj}) = 0$, like in (2). ■

As a result of Proposition 1, model (5) can be reformulated as follows, wherein we have removed the absolute values in the objective function by applying the change of variables normally used to that end

$$
\begin{array}{llll}
\text{Min} & \displaystyle\sum_{j\in J_g}\sum_{r=1}^{s}\frac{\left(p_{rj}^+ + p_{rj}^-\right)}{Q_j\omega_r} & & \\
\text{s.t.:} & & & \\
 & \text{constraints (5.1)-(5.4)} & & \\
 & Q_j\omega_r I_{rj}^1 + Q_j\omega_r\left(1 - \frac{s_{rj}}{y_{rj}}\right)I_{rj}^2 + \left(p_{rj}^+ - p_{rj}^-\right) = p_{rj}\left(s_{rj}^g\right) & j \in J_g \quad r = 1,\dots,s & (7.1) \\
 & y_{rj}I_{rj}^3 - M_1 I_{rj}^1 \le s_{rj} \le y_{rj}I_{rj}^2 + M_2 I_{rj}^3 & j \in J_g \quad r = 1,\dots,s & (7.2) \\
 & I_{rj}^1 + I_{rj}^2 + I_{rj}^3 = 1 & j \in J_g \quad r = 1,\dots,s & (7.3) \\
 & \lambda_k = \displaystyle\sum_{j\in J_g}\lambda_{kj} & k \in E & (7.4) \\
 & p_{rj}^+, p_{rj}^- \ge 0 & j \in J_g \quad r = 1,\dots,s & \\
 & d_k \ge 0, \lambda_{kj} \ge 0 & k \in E \quad j \in J_g & \\
 & \{\lambda_k, d_k\} \quad \text{SOS1} & k \in E & \\
 & I_{rj}^1, I_{rj}^2, I_{rj}^3 \in \{0,1\} & j \in J_g \quad r = 1,\dots,s & \\
 & u_0, s_{rj} \quad \text{free} & j \in J_g \quad r = 1,\dots,s & (7)
\end{array}
$$

The following proposition allows us to find lower bounds for $M_1$ and $M_2$

**Proposition 2.** For every DMU$_j$, $j \in J_g$, and each indicator $y_r$, r=1,…,s, the following statement holds

$$-y_{rj} \le s_{rj} \le \max_{k\in E}\{y_{rk}\}$$

Proof. From (5.1), we have $y_{rj} + s_{rj} = \sum_{k \in E} \lambda_{kj} y_{rk} \geq 0$. Then, $s_{rj} \geq -y_{rj}$. On the other hand, $s_{rj} \leq y_{rj} + s_{rj} = \sum_{k \in E} \lambda_{kj} y_{rk} \leq \max_{k \in E} \{y_{rk}\}$, because of (5.2). ■

**Corollary 1.** It suffices to set $M_1$ at $y_{rj}$ in the constraint (7.2) corresponding to indicator $y_r$ of $DMU_j$ and $M_2$ at $\max_{k \in E} \{y_{rk}\}$ in the constraints involving indicator $y_r$.

Model (7) is still a non-linear problem because of the products of continuous and binary variables $s_{rj} I_{rj}^2$, $j \in J_g$, $r = 1,...,s$. Using the change of variables $z_{rj} = s_{rj} I_{rj}^2$, $j \in J_g$, $r = 1,...,s$, an equivalent mixed-integer linear programming problem (MILP) is formulated in the next proposition, wherein we have also considered the lower bounds for $M_1$ and $M_2$ which were found in corollary 1.

**Proposition 3.** Model (7) is equivalent to the following MILP problem

$$
\begin{array}{llll}
\text{Min} & \displaystyle\sum_{j \in J_g} \sum_{r=1}^{s} \frac{\left(p_{rj}^{+} + p_{rj}^{-}\right)}{Q_j \omega_r} & & \\
\text{s.t.:} & & & \\
 & \text{constraints (5.1)-(5.4)} & & \\
 & Q_j \omega_r I_{rj}^1 + Q_j \omega_r I_{rj}^2 - Q_j \omega_r \dfrac{z_{rj}}{y_{rj}} + \left(p_{rj}^{+} - p_{rj}^{-}\right) = p_{rj}\left(s_{rj}^{g}\right) & j \in J_g \quad r = 1,...,s & (8.1) \\
 & -y_{rj} I_{rj}^2 \leq z_{rj} \leq I_{rj}^2 \max_{k \in E} \{y_{rk}\} & j \in J_g \quad r = 1,...,s & (8.2) \\
 & -y_{rj}(1 - I_{rj}^2) \leq s_{rj} - z_{rj} \leq (1 - I_{rj}^2) \max_{k \in E} \{y_{rk}\} & j \in J_g \quad r = 1,...,s & (8.3) \\
 & y_{rj}(I_{rj}^3 - I_{rj}^1) \leq s_{rj} \leq y_{rj} I_{rj}^2 + I_{rj}^3 \max_{k \in E} \{y_{rk}\} & j \in J_g \quad r = 1,...,s & (8.4) \\
 & I_{rj}^1 + I_{rj}^2 + I_{rj}^3 = 1 & j \in J_g \quad r = 1,...,s & (8.5) \\
 & \lambda_k = \displaystyle\sum_{j \in J_g} \lambda_{kj} & k \in E & (8.6) \\
 & p_{rj}^{+}, p_{rj}^{-} \geq 0 & j \in J_g \quad r = 1,...,s & \\
 & d_k \geq 0, \lambda_{kj} \geq 0 & k \in E \quad j \in J_g & \\
 & \{\lambda_k, d_k\} \quad \text{SOS1} & k \in E & \\
 & I_{rj}^1, I_{rj}^2, I_{rj}^3 \in \{0,1\} & j \in J_g \quad r = 1,...,s & (8) \\
 & u_0, s_{rj}, z_{rj} \quad \text{free} & j \in J_g \quad r = 1,...,s &
\end{array}
$$

Proof. It suffices to prove that (7) and (8) have the same constraints in the different scenarios resulting from the specification of the binary variables.

Suppose that $I^1_{rj} = 1$ for some r and some j. Then, the corresponding constraints for these r and j in (7.1) and (7.2) become $Q_j\omega_r + \left(p^+_{rj} - p^-_{rj}\right) = p_{rj}\left(s^g_{rj}\right)$ and $-y_{rj} \le s_{rj} \le 0$, respectively, while in (8.1)-(8.4) those constraints would be, respectively, $Q_j\omega_r + \left(p^+_{rj} - p^-_{rj}\right) = p_{rj}\left(s^g_{rj}\right)$, $z_{rj} = 0$, $-y_{rj} \le s_{rj} \le \max_{k\in E}\left\{y_{rk}\right\}$ and $-y_{rj} \le s_{rj} \le 0$. That is, both sets of constraints coincide.

If $I^2_{rj} = 1$, then in (7) we have $Q_j\omega_r - \frac{s_{rj}}{y_{rj}} Q_j\omega_r + \left(p^+_{rj} - p^-_{rj}\right) = p_{rj}\left(s^g_{rj}\right)$ and $0 \le s_{rj} \le y_{rj}$, whereas in (8) (8.3) implies $s_{rj} = z_{rj}$, so (8.1) becomes $Q_j\omega_r - \frac{s_{rj}}{y_{rj}} Q_j\omega_r + \left(p^+_{rj} - p^-_{rj}\right) = p_{rj}\left(s^g_{rj}\right)$, (8.2) is $-y_{rj} \le s_{rj} \le \max_{k\in E}\left\{y_{rk}\right\}$ and (8.4) is $0 \le s_{rj} \le y_{rj}$.

Finally, if $I^3_{rj} = 1$, (7.1) and (7.2) become $p^+_{rj} - p^-_{rj} = p_{rj}\left(s^g_{rj}\right)$ and $y_{rj} \le s_{rj} \le \max_{k\in E}\left\{y_{rk}\right\}$, respectively, while in (8) (8.1)-(8.4) are, respectively, $p^+_{rj} - p^-_{rj} = P_{rj}\left(s^g_{rj}\right)$, $z_{rj} = 0$, $-y_{rj} \le s_{rj} \le \max_{k\in E}\left\{y_{rk}\right\}$ and $y_{rj} \le s_{rj} \le \max_{k\in E}\left\{y_{rk}\right\}$. ■

Eventually, we can obtain the optimal solutions of model (4) by solving problem (8). Note that (8) is a mixed integer linear problem (MILP) in which the number of binary variables and SOS variables depends on the number of indicators considered and that of efficient DMUs. In the DEA applications that we usually find in practice (8) becomes a problem that can be easily solved by using conventional software.

## 4. Illustrative example

For purposes of illustration only, in this section we apply the proposed approach in an example concerning the financing of public Spanish universities.

The public Spanish universities may be seen as a set of homogeneous DMUs that undertake similar activities and produce comparable results in respect of their missions (teaching, research,…), so that they can be compared in terms of a common set of outputs which can be defined for the analysis (see Dyson et al. (2001) for a discussion on homogeneity assumptions about the units under assessment). Public Spanish universities are not specialized either in teaching or research. They all come under the same set of general regulations. Nevertheless, in Spain, many of the competences related to education are transferred to the regional (autonomous) governments. This means that the management in these Higher Education Institutions (HEI) is determined to some extent by the policies established by these local governments. These include the offer of study programs, the fixing of university fees, the policy for the promotion of research, the setting of goal levels for quality indicators of teaching and research, the financing models and even the regulation of supplements for staff salaries.

Concerning the financing of universities, the models used in the different regions are distinguished mainly by their components of base financing and other issues of a more strategic nature. Base financing includes government grants and own resources of the universities, like tuition fees and research incomes, and is aimed at covering ordinary running costs, delivery of educational services and investments. However, universities can often have access to extra financing by means of the implementation of pay-for-performance incentive plans, which are typically linked to the achievement of goals. Thus, policy makers and university managers of HEI, through this strategic component of financing, seek to offer ways that may help universities to achieve the desired aims of efficiency, effectiveness and quality.

Incentive plans at public Spanish universities typically pay for performance in areas like teaching or research, which is evaluated in terms of some index indicators. Payment is made on the basis of the degree of achievement of goals. As for the setting of goal levels for each university, these often result from an agreement (a programme-contract) between the corresponding regional government and their own managers of that university, in which stakeholders take into consideration their objectives and policies, the knowledge of prior period performance, etc…

This situation of HEI in Spain, where, as said before, many of the competences in education are transferred to the regional governments, suggests that the universities of a given region are experiencing similar circumstances for the development of the activities in the context of their different missions. Therefore, it would be desirable to treat them uniformly in the benchmarking analysis that is carried out for the setting of targets. For this reason, we use the approach proposed in this paper by considering as groups the universities in each of the regions, so that targets result from their benchmarking against a common reference set (region-specific) of efficient universities, as a consequence of not allowing for individual circumstances within the autonomous regions. Note, however, that this DEA-based approach allows for regional circumstances, which would be justified by the possible variations in their goals and policies.

The analysis carried out here is concerned with the overall performance of the universities in the sense that indicators related to the two main areas of activity, teaching and research, are considered. In order to do so, we have selected both indicators such as rates of progress, graduation and retention, and others related to publications (production and quality) and research income. Specifically, the following indicators have been considered for the analysis:

- PROGRESS RATE: Ratio of the number of passed credits[2] corresponding to all the students enrolled in 2013-14, to the total enrolled credits in that academic year (as a percentage).
- GRADUATION RATE: Percentage of students that completed the programme of studies in 2013-14, in relation to number of students enrolled for the first time in 2010-11.
- RETENTION RATE: This is computed as 100 minus the drop out rate, in order to be treated as an output, i.e., a "the more the better" variable. The drop out rate is the percentage of students that abandon in 2013-14 the subject in which they enrolled for the first time in 2011-12.
- ARTICLES: Ratio of the total number of articles published in the period 2009-13, to the academic staff in the academic year 2013-14 (with full-time equivalence).
- %Q1: Percentage of articles published in journals classified in the 1st quartile of the corresponding category.
- RESEARCH INCOME: Ratio of total income from research projects (grants) in competitive calls and contracts with companies and administrations in 2013 to academic staff in the academic year 2013-14 (with full-time equivalence).

Data for these variables have been taken from the reports by the Conference of Rectors of the Spanish Universities (CRUE), and by the Foundation FCyD[3]. Table 3 records the universities considered in the analysis, classified by regions (the selection of the sample has been determined by the availability of data[4]).

Table 3. Universities by regions

The DEA analysis of the entire sample reveals 9 efficient universities: UAB, UBA, UDL, UPF, URI, UVA, UAM, UMH and UPVA. The universities in each of the regions are benchmarked against reference sets associated with some of the facets of the Pareto frontier formed by these 9 universities. We report the results corresponding to the benchmarking analyses of three regions as representative cases: Cataluña, Andalucía and Comunidad Valenciana. These are recorded in Tables

[2] Credit is the unit of measurement of the academic load of the subject of a programme.
[3] It could be translated as Foundation for the Knowledge and Development.
[4] The rate of graduation for UVEG has been estimated.

4, 5 and 6, respectively, together with actual data. Each of these three tables contains information regarding both the DEA benchmarking and the incentive plans. We suppose that, in the academic year 2012-13, goals were set for each university as specified in the rows "goal" of these tables. The last column of the tables records the initially available amount of money for each university, in the incentive plan (in the first row corresponding to each university). This information has also been taken from the CRUE report and represents approximately 10% of the current transfers coming from the regional governments in 2013. We also report the money that the universities receive, on the basis of the evaluation of their performance relative to the DEA targets and the goals, and the payment rates (expressed in terms of percentages of the available amount). As for the benchmarking, the tables record, for each university, its actual data and the DEA targets corresponding to each indicator (as obtained from model (8)), aside from the goals.

Table 4. DEA target setting and goals (Cataluña)

Table 5. DEA target setting and goals (Andalucía)

Table 6. DEA target setting and goals (Comunidad Valenciana)

Looking at the payment rates allows us to identify universities which still have room for improvement after the implementation of the incentive plans. In this respect, it should be pointed out that the payments resulting from the DEA targets are very similar to those relative to the goals in all of the universities. This would validate the payments relative to the goals initially established, because the amounts of money given to the universities are practically the same as those that would be derived from an evaluation of performance relative to targets which are attainable and efficient. Nevertheless, DEA targets and goals differ for some indicators in some universities. DEA targets are sometimes higher and sometimes lower than the goals. Only in the case of the RETENTION rate, are DEA targets always more demanding. We note, in addition, that goals (the goal bundles) for UAB, UPC, UPF and UPVA, represent plans of performance that are not achievable, because they do not belong to the attainable set AS. For example, in the case of UPC, the target set by the model is lower than the goal for ARTICLES, while at the same time that of INCOMES has been set at a more demanding level (in addition to making some small adjustments in other indicators), in order to set an attainable bundle of DEA targets for the evaluation of this university. For all these reasons, we carry out the benchmarking by using the DEA targets.

Tables show a better performance in the case of universities of Cataluña. On average, the universities in Cataluña have a payment rate relative to DEA targets of 90.28%, while in Andalucía this rate is of 80.78%, even though the levels of the targets set in that region are lower. In Cataluña,

whose universities are benchmarked against UPF, UVA, UAM, UMH and UPVA, we highlight the case of UPF, whose payment rate is of 98.4%. This means that this university gets practically the entire amount of money available in the incentive plan. Its actual performance exceeds that arising from the DEA targets in almost all the indicators, except in the RETENTION rate wherein this university should improve. In contrast, UPC shows a poorer performance (its payment rate is 75.7%). DEA targets reveal that this university performs well in terms of the numbers of ARTICLES published, but it should improve significantly the quality of its publications (%Q1); in INCOMES it has achieved the objective. As for teaching performance, UPC has an important weakness in the GRADUATION of students (the actual rate, 13.21%, is far below the target of 36.88%), while some small improvements are also needed in the other rates.

In Andalucía, where the reference set is formed by UPF, URI, UVA, UAM, UMH and UPVA, UAL and UGR are the universities that have performed better, albeit they have some important weaknesses: UAL in INCOMES and UGR in the rate of GRADUATION. In contrast, UCA and UMA, whose payment rates are 70.5% and 73.6%, respectively, still have much room for improvement both in teaching and research activities. Finally, in Comunidad Valenciana, UPVA and UMH have been the best performers, while UA and UJI should improve the production of ARTICLES and INCOMES aside from the rate of GRADUATION.

## 5. Conclusions

The DEA literature dealing with the evaluation of performance in presence of management goals is scarce. This paper makes a contribution in this area through an approach that provides a benchmarking framework based on appropriate targets for the monitoring of improvements of DMUs that follow the implementation of pay-for-performance incentive plans. In this context, the goals that have been set may sometimes be unachievable or, on the contrary, may not be ambitious enough (sometimes because of participative goal setting). We show how DEA can be extended through a process of benchmarking adjusted to the goals that allows us to set targets that are attainable and represent best practices. The setting of targets also takes into account the policy of improvements that is pursued through the implementation of incentive plans, because we seek targets that respect, as much as possible, the payments of incentives relative to the goals set.

The proposed methodology should be of special interest in areas where pay-for-performance plans are usually implemented, like Health Care. In this area, DEA has already been used to create a composite measure of quality within the framework of a pay-for-performance program (see Shwartz et al., 2016). Our approach would provide an alternative for performance evaluation through the setting of targets that represent best practice performances. Higher Education is another potential area of application of the methodology, as has been illustrated here with an example in which we evaluate performance of public Spanish universities in the context of their financing. Nigsch and Schenker-Wicki (2015) raise the need to use efficiency analysis for incentive-setting and the reallocation of funds, in a discussion on the use of DEA for analyzing the performance of universities.

Regarding future directions, a couple of lines of research can be envisaged. On one hand, we should develop a dynamic DEA framework, where targets and goals are adjusted with incentives and penalties over a several time periods. On the other, new DEA benchmarking models for performance evaluation should be formulated within a more general framework, wherein management goals are set, but without having a monetary incentive for their achievement (in line with those proposed in Stewart (2010)).

## Acknowledgments

This research has been supported through Grant MTM2016-76530-R (AEI/FEDER, UE).

Table 1. Data of numerical example

| | Period t-1 | | | Period t | |
|---|---|---|---|---|---|
| DMU | goal y1 | goal y2 | Total Incentives | y1 | y2 |
| A | 3 | 7 | 25 | 1 | 7 |
| B | 5 | 4 | 30 | 6 | 5 |
| C | 8 | 4 | 20 | 9 | 1 |
| D | 4 | 7 | 20 | 3 | 4 |
| E | 6 | 3 | 25 | 5 | 2 |
| F | 2 | 6.6 | 20 | 2 | 5 |

Table 2. Results of numerical example

| | | Data | | Payments of incentives | | | |
|---|---|---|---|---|---|---|---|
| DMU | | y1 | y2 | | y1 | y2 | Total |
| A | Actual | 1 | 7 | Available | 12.5 | 12.5 | 25 |
| | Targets | 2 | 6.6 | Targets | 0 | 12.5 | 12.5 |
| | Goals | 3 | 7 | Goals | 0 | 12.5 | 12.5 |
| B | Actual | 6 | 5 | Available | 15 | 15 | 30 |
| | Targets | 6 | 5 | Targets | 15 | 15 | 30 |
| | Goals | 5 | 4 | Goals | 15 | 15 | 30 |
| C | Actual | 9 | 1 | Available | 10 | 10 | 20 |
| | Targets | 8.25 | 2 | Targets | 10 | 0 | 10 |
| | Goals | 8 | 4 | Goals | 10 | 0 | 10 |
| D | Actual | 3 | 4 | Available | 10 | 10 | 20 |
| | Targets | 4 | 5.8 | Targets | 6.67 | 5.5 | 12.17 |
| | Goals | 4 | 7 | Goals | 6.67 | 2.5 | 9.17 |
| E | Actual | 5 | 2 | Available | 12.5 | 12.5 | 25 |
| | Targets | 7.5 | 3 | Targets | 6.25 | 6.25 | 12.5 |
| | Goals | 6 | 3 | Goals | 10 | 6.25 | 16.25 |
| F | Actual | 2 | 5 | Available | 10 | 10 | 20 |
| | Targets | 2 | 6.6 | Targets | 10 | 6.8 | 16.8 |
| | Goals | 2 | 6.6 | Goals | 10 | 6.8 | 16.8 |

Table 3. Universities by regions

| | |
|---|---|
| **ANDALUCÍA** | **CANARIAS** |
| U. de ALMERÍA (UAL) | U. de LA LAGUNA (ULL) |
| U. de CÁDIZ (UCA) | U. de LAS PALMAS de G. C. (ULPGC) |
| U. de GRANADA (UGR) | |
| U. de JAÉN (UJA) | **ARAGÓN** |
| U. de MÁLAGA (UMA) | U. de ZARAGOZA (UZA) |
| U. de SEVILLA (USE) | |
| | **ASTURIAS** |
| **CATALUÑA** | U. de OVIEDO (UOV) |
| U. AUTÓNOMA de BARCELONA (UAB) | |
| U. de BARCELONA (UBA) | **CANTABRIA** |
| U. de GIRONA (UDG) | U. de CANTABRIA (UCN) |
| U. de LLEIDA (UDL) | |
| U. POLITÉCNICA DE CATALUÑA (UPC) | **CASTILLA LA MANCHA** |
| U. POMPEU FABRA (UPF) | U. de CASTILLA-LA MANCHA (UCLM) |
| U. ROVIRA I VIRGILI (URV) | |
| | **EXTREMADURA** |
| **CASTILLA Y LEÓN** | U. de EXTREMADURA (UEX) |
| U. de BURGOS (UBU) | |
| U. de LEÓN (ULE) | **ISLAS BALEARES** |
| U. de VALLADOLID (UVA) | U. de las ISLAS BALEARES (UIB) |
| | |
| **COMUNIDAD VALENCIANA** | **MURCIA** |
| U. de ALICANTE (UA) | U. de MURCIA (UMU) |
| U. JAUME I de CASTELLÓN (UJCS) | U. POLITÉCNICA DE CARTAGENA (UPCT) |
| U. MIGUEL HERNÁNDEZ de ELCHE (UMH) | |
| U. POLITÉCNICA DE VALENCIA (UPVA) | **NAVARRA** |
| U. de VALENCIA (ESTUDI GENERAL) (UVEG) | U. PÚBLICA DE NAVARRA (UPN) |
| | |
| **MADRID** | **PAIS VASCO** |
| U. de ALCALÁ de HENARES (UAH) | U. del PAÍS VASCO (UPV) |
| U. AUTÓNOMA de MADRID (UAM) | |
| | **LA RIOJA** |
| **GALICIA** | U. de LA RIOJA (URI) |
| U. de LA CORUÑA (ULC) | |
| U. de VIGO (UVI) | |

Table 4. DEA target setting and goals (Cataluña)

| | | ARTICLES | %Q1 | INCOMES | PROGRESS | GRADUATION | RETENTION | Payments |
|---|---|---|---|---|---|---|---|---|
| UAB | actual | 7.14 | 59.48 | 14256.16 | 84.39 | 36.08 | 81.95 | 15,000,000.00 |
| | target | 5.65 | 59.90 | 15000 | 85 | 40 | 85.45 | 14,455,127.27 (96.4%) |
| | goal* | 6 | 60 | 15000 | 85 | 40 | 85 | 14,464,537.02 (96.4%) |
| UBA | actual | 5.59 | 62.29 | 8383.88 | 83.98 | 29.79 | 81.55 | 20,000,000.00 |
| | target | 4.06 | 53.03 | 10000 | 82.69 | 40 | 88.30 | 17,938,771.89 (89.7%) |
| | goal | 5 | 60 | 10000 | 85 | 40 | 85 | 18,033,392.57 (90.2%) |
| UDG | actual | 4.15 | 52.27 | 11815.68 | 83.55 | 29.39 | 77.24 | 5,000,000.00 |
| | target | 4.15 | 52.14 | 11815.68 | 82.74 | 40 | 87.86 | 4,584,635.52 (91.7%) |
| | goal | 4 | 50 | 10000 | 85 | 40 | 85 | 4,601,029.47 (92.0%) |
| UDL | actual | 3.20 | 58.10 | 6492.60 | 81.91 | 43.03 | 85.40 | 4,000,000.00 |
| | target | 4 | 52.66 | 10000 | 82.62 | 38.94 | 88.41 | 3,442,643.38 (86.1%) |
| | goal | 4 | 60 | 10000 | 85 | 40 | 85 | 3,424,922.83 (85.6%) |
| UPC | actual | 6.68 | 39.45 | 17847.51 | 76.05 | 13.21 | 78.94 | 13,000,000.00 |
| | target | 4.55 | 50 | 17847.51 | 83.42 | 36.88 | 86.32 | 9,841,404.20 (75.7%) |
| | goal* | 6 | 50 | 15000 | 85 | 40 | 85 | 9,832,693.22 (75.6%) |
| UPF | actual | 8.12 | 59.8 | 35555.19 | 89.95 | 51.66 | 77.95 | 5,000,000.00 |
| | target | 5.69 | 60 | 15702.07 | 86.58 | 40.49 | 85 | 4,921,862.56 (98.4%) |
| | goal* | 6 | 60 | 25000 | 85 | 40 | 85 | 4,921,862.56 (98.4%) |
| URV | actual | 4.76 | 57.83 | 12250.08 | 83.14 | 41.21 | 80.88 | 6,000,000.00 |
| | target | 5.00 | 55.61 | 15000 | 84.83 | 39.06 | 86.13 | 5,639,714.23 (94.0%) |
| | goal | 5 | 60 | 15000 | 85 | 40 | 85 | 5,614,153.06 (93.6%) |

* Goals for this university represent a plan which is outside AS

Table 5. DEA target setting and goals (Andalucía)

| | | ARTICLES | %Q1 | INCOMES | PROGRESS | GRADUATION | RETENTION | Payments |
|---|---|---|---|---|---|---|---|---|
| UAL | actual | 3.68 | 44.50 | 11903.49 | 77.70 | 35.08 | 80.44 | 6,000,000.00 |
| | target | 4 | 50 | 15000 | 80.43 | 32.18 | 87.81 | 5,402,694.12 (90.0%) |
| | goal | 4 | 50 | 15000 | 80 | 35 | 85 | 5,443,169.03 (90.7%) |
| UCA | actual | 1.83 | 48.23 | 8961.71 | 77.64 | 17.95 | 80.89 | 11,000,000.00 |
| | target | 3 | 50 | 6170.32 | 79.24 | 35 | 90.60 | 7,759,069.88 (70.5%) |
| | goal | 3 | 50 | 10000 | 80 | 35 | 85 | 7,655,603.61 (69.6%) |
| UGR | actual | 4.05 | 52.37 | 15241.30 | 77.24 | 22.29 | 80.98 | 17,000,000.00 |
| | target | 4.05 | 55 | 8587.71 | 80.17 | 35 | 89.00 | 14,853,727.62 (87.4%) |
| | goal | 4 | 55 | 15000 | 80 | 35 | 85 | 14,999,711.36 (88.2%) |
| UJA | actual | 3.33 | 48.12 | 9157.57 | 77.11 | 22.94 | 78.84 | 7,000,000.00 |
| | target | 4 | 51.76 | 10000 | 78.53 | 35 | 89.20 | 5,782,554.31 (82.6%) |
| | goal | 4 | 50 | 10000 | 80 | 35 | 85 | 5,865,140.54 (83.8%) |
| UMA | actual | 2.67 | 44.18 | 9622.92 | 75.05 | 27.80 | 82.33 | 11,000,000.00 |
| | target | 4 | 49.46 | 15000 | 81.20 | 35 | 87.67 | 8,095,434.16 (73.6%) |
| | goal | 4 | 50 | 15000 | 80 | 35 | 85 | 8,161,669.27 (74.2%) |
| USE | actual | 3.02 | 52.71 | 13723.68 | 72.24 | 22.24 | 82.49 | 18,000,000.00 |
| | target | 4 | 47.96 | 15000 | 80.11 | 35 | 87.99 | 14,497,122.07 (80.5%) |
| | goal | 4 | 55 | 15000 | 80 | 35 | 85 | 14,480,276.99 (80.4%) |

Table 6. DEA target setting and goals (Comunidad Valenciana)

| | | ARTICLES | %Q1 | INCOMES | PROGRESS | GRADUATION | RETENTION | Payments |
|---|---|---|---|---|---|---|---|---|
| UA | actual | 2.76 | 48.39 | 6868.54 | 81.02 | 30.38 | 89.92 | 12,000,000.00 |
| | target | 4 | 52.60 | 10000 | 82.52 | 40 | 88.38 | 9,346,396.51 (77.9%) |
| | goal | 4 | 50 | 10000 | 85 | 40 | 85 | 9,392,607.97 (78.3%) |
| UJI | actual | 2.97 | 48.88 | 7810.98 | 77.91 | 28.42 | 80.70 | 7,000,000.00 |
| | target | 4 | 52.60 | 10000 | 82.52 | 40 | 88.38 | 5,526,258.16 (78.9%) |
| | goal | 4 | 50 | 10000 | 85 | 40 | 85 | 5,600,097.60 (80.0%) |
| UMH | actual | 4.36 | 54.34 | 8703.58 | 74.56 | 31.51 | 89.86 | 6,000,000.00 |
| | target | 4.24 | 54.34 | 10000 | 83.21 | 39.80 | 88.05 | 5,471,761.85 (91.2%) |
| | goal | 4 | 50 | 10000 | 85 | 40 | 85 | 5,441,556.96 (90.7%) |
| UPVA | actual | 4.97 | 44.96 | 27320.81 | 83.72 | 33.82 | 84.18 | 21,000,000.00 |
| | target | 5 | 50 | 21387.06 | 83.96 | 40 | 85.02 | 19,901,068.16 (94.8%) |
| | goal* | 5 | 50 | 25000 | 85 | 40 | 85 | 19,858,344.15 (94.6%) |
| UVEG | actual | 3.02 | 52.71 | 13723.68 | 72.24 | 22.24 | 82.49 | 24,000,000.00 |
| | target | 5 | 55.56 | 15000 | 84.74 | 40 | 86.11 | 21,368,626.46 (89.0%) |
| | goal | 5 | 50 | 15000 | 85 | 40 | 85 | 21,357,035.96 (89.0%) |

* Goals for this university represent a plan which is outside AS

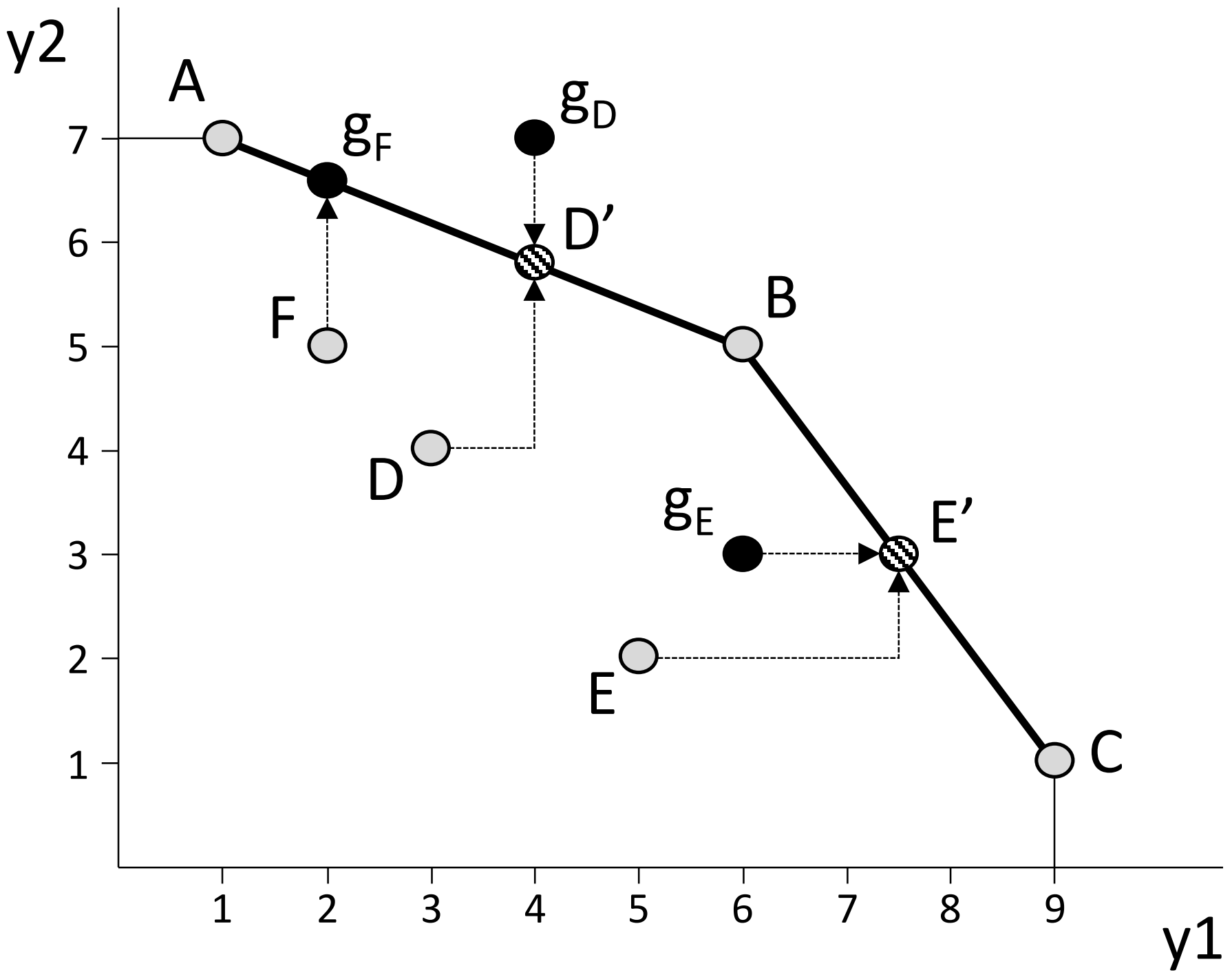

Figure 2. Target setting for DMUs D, E and F

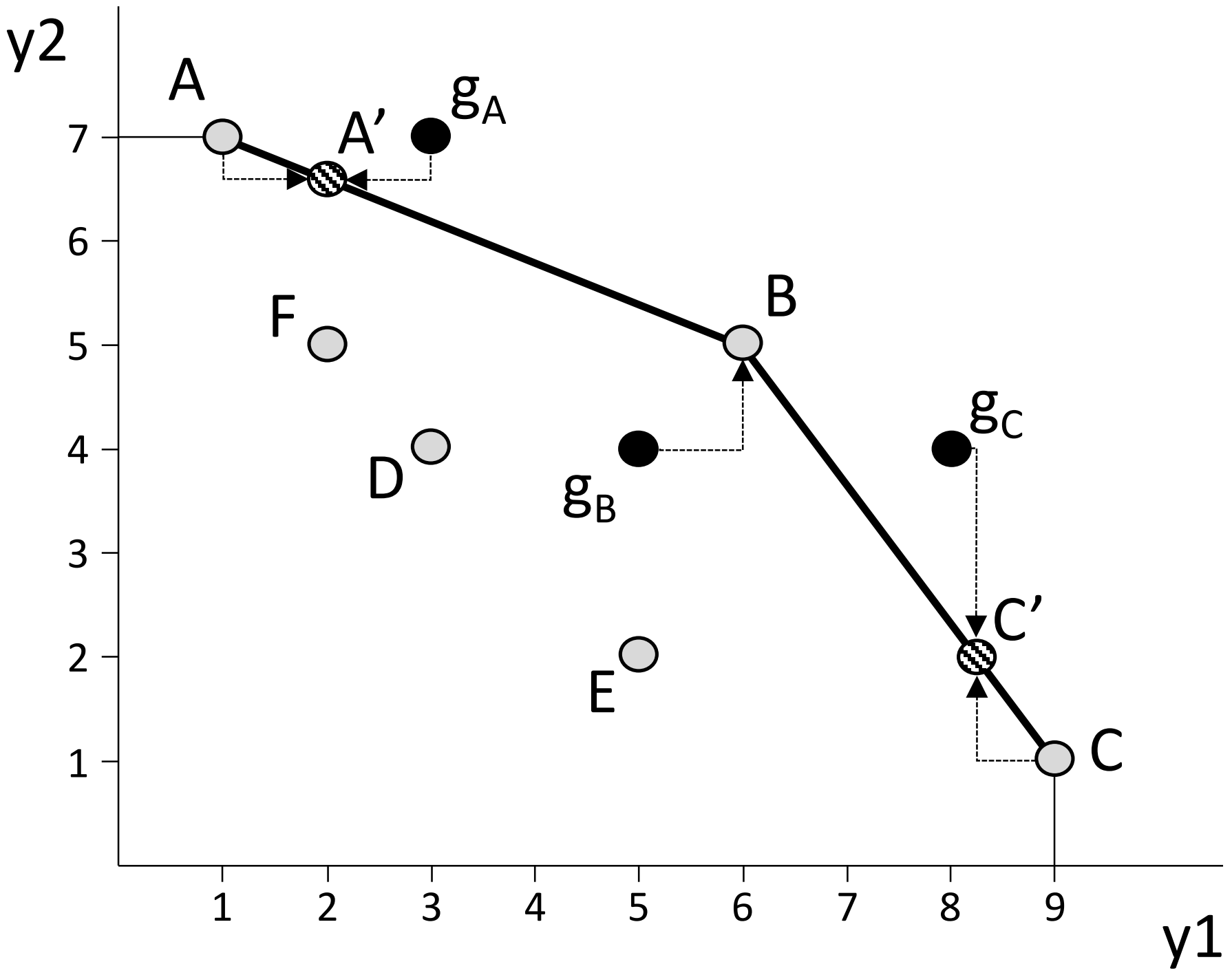

Figure 3. Target setting for DMUs A, B and C